\documentclass[11pt,twoside]{article}
\usepackage{latexsym,amsmath}
\usepackage{indentfirst}
\usepackage{graphicx}

\numberwithin{equation}{section}
\newtheorem{Prop}{\bf Proposition}[section]
\newtheorem{Cor}{\bf Corollary}[section]

\newtheorem{Rem}{\bf Remark}[section]
\newtheorem{Ex}{\bf Example}[section]

\begin{document}
\def \b{\Box}
\begin{center}
{\Large {\bf On the 3-dimensional fractional-order Toda lattice with two controls}}
\end{center}

\begin{center}
{\bf Mihai IVAN}\\
\end{center}

\setcounter{page}{1}
\pagestyle{myheadings}

{\small {\bf Abstract}. The main purpose of this paper is  to study the fractional-order system  with Caputo derivative associated to 3-dimensional Toda lattice with two controls. For this fractional system we investigate the existence and uniqueness of solution of initial value problem, asymptotic stability of its equilibrium states,  stabilization problem using appropriate controls and numerical integration via the fractional Euler method.}
{}

\section{Introduction}

Fractional calculus is a interesting field of mathematics based on a generalization of integer derivatives to fractional-order. The study of this field is very important and attractive because several phenomena have been proven to be better described by fractional derivatives that take into account not only the local properties,  but also global correlations of dynamical systems \cite{agra, nabu, mati}. The fractional calculus has deep and natural connections with many fields of science and engineering \cite{diet, kilb, podl}.

In the last three decades, one increasing attention has been paid to the study of the dynamic behaviors (in particular, the chaotic behavior) of some classical differential systems, as well as some fractional-order differential systems. For example, the fractional models played an important role in  applied mathematics \cite{enpo, ivan, ghiv, igim, opmi, migi}, mathematical physics \cite{ivmi, mihi}, applied physics \cite{mi22, miva24, mivan24, puca, sczc},  study of biological systems  \cite{ahme,lima, pagi}, chaos synchronization, secure communications \cite{bhda,igmp, kich} and so on.

The  topological groupoids,  Lie groupoids,  Lie algebroids and Leibniz algebroids have proven to be powerful tools for geometric formulation of the Hamiltonian mechanics \cite{demi, giop,gimo, miv03, miva05}. Also, they have been used  in the investigation of many fractional dynamical systems \cite{mihi, migo}.

A new generalization of the notion of groupoid \cite{ivan02, ivan24, miva03,miva13} is defined by Mihai Ivan in paper \cite{miva23a}, called almost groupoid. This algebraic structure can be used to define the concepts of almost topological groupoid and almost Lie groupoid.

This paper is structured as follows. In Section 2, we investigate the fractional differential systems associated to 3-dimensional Toda lattice with two controls $ (2.4)~$ in terms of fractional Caputo derivatives. This new fractional system $ (2.4) $ is called the 3-dimensional fractional-order Toda lattice with two controls. The existence and uniqueness of solution of initial value problem for the fractional system $ (2.4) $ are proved. In Propositions $ 2.4 $ and $ 2.5, $ the asymptotic stability for the equilibrium states of the fractional system $ (2.4) $ is discussed. In Section $ 3,$ for the stability controlability problem we associate the fractional system  $ (3.3),$ called the controlled fractional-order Toda lattice associated to $ (2.4) $ at an equilibrium state.  In Propositions $~3.1-3.2~$ and Corollaries $~3.1-3.3~$ are established sufficient conditions on parameters  $~a, b, c_{1}, c_{2}, c_{3}~$  to control the chaos in the fractional system $~(3.3).$ The numerical integration of the model $ (3.3)~$ based on fractional Euler method, is presented in Section 4.\\[-0.5cm]

\section{ Stability analysis of the $ 3-$dimensional fractional-order Toda lattice with two controls}

We recall the Caputo definition of fractional derivatives, which is often used in concrete applications. Let $ f\in C^{\infty}(
\textbf{R}) $ and $ q \in \textbf{R}, q > 0. $ The $ q-$order Caputo differential operator
\cite{diet}, is described by $~D_{t}^{q}f(t) = I^{m - q}f^{(m)}(t), ~q > 0,~$ where $~f^{(m)}(t)$  represents the $
m-$order derivative of the function $ f,~m \in \textbf {N}^{\ast}$
is an integer such that $ m-1 \leq q \leq m $ and $ I^{q} $ is the
$ q-$order Riemann-Liouville integral operator, which
is expressed by $~I_{t}^{q}f(t)
=\displaystyle\frac{1}{\Gamma(q)}\int_{0}^{t}{(t-s)^{q
-1}}f(s)ds,~q > 0,$ where $~\Gamma $ is the Euler Gamma  function.
 If $ q =1$, then $ D_{t}^{q}f(t) = df/dt.~$  In this paper we suppose that $ q \in (0,1].$

The Toda lattice is a well-known model for a one-dimensional crystal in solid state physics \cite{toda}.

In the Euclidean space $ {\bf R}^{2n-1}~$ with local coordinates $\{ x^{1}, x^{2},\ldots, x^{n}, y^{1}, y^{2}, \ldots, y^{n-1}\} $ the {\it $~n-$dimensional Toda lattice} dynamics is described by the following set of differential equations \cite{dami, pudo}:\\[-0.2cm]
\begin{equation}
\dot{x}^{i}(t)  =  2[(y^{i})^{2}(t)-(y^{i-1})^{2}(t)], ~~~
\dot{y}^{j}(t)  =  y^{j}(t)( x^{j+1}(t)-x^{j}(t)),~~~\label{(2.1)}
\end{equation}
where $~y^{0}(t) = y^{n}(t)=0,~x^{i} $ for $ i=\overline{1,n}, $ and $~y^{j} $ for $ j=\overline{1,n-1}~$ are state variables, $~\dot{x}^{i}(t)= dx^{i}(t)/dt,~\dot{y}^{j}(t)= dy^{j}(t)/dt~$ and $ t $ is the time.

The {\it n- dimensional fractional-order Toda lattice} associated to dynamics $~(2.1)~$ is defined by the following
set of fractional differential equations:\\[-0.2cm]
\begin{equation}
\left\{ \begin{array} {lll}
 D_{t}^{q}{x}^{i}(t) & = &  2[(y^{i})^{2}(t)-(y^{i-1})^{2}(t)],~~~~~~i=\overline{1,n}  \\[0.1cm]
 D_{t}^{q}{y}^{j}(t) & = & y^{j}(t)( x^{j+1}(t)-x^{j}(t))~~~~~~~~j=\overline{1,n-1},~~~~~  q \in (0,1),\\[0.1cm]
  y^{0}(t) = 0, & & y^{n}(t)=0. \label{(2.2)}
  \end{array}\right.
\end{equation}
Let us we consider  the dynamics of the n-dimensional fractional-order Toda lattice for $~n=3 $ with two particular controls around axes $~Ox^{1}~$ and $~Oy^{2}.~$ This fractional model is described  by:\\[-0.2cm]
\begin{equation}
\left\{ \begin{array} {lll}
D_{t}^{q}{x}^{1}(t) & = & 2 (y^{1})^{2}(t) + a x^{1}(t),\\[0.1cm]
D_{t}^{q}{x}^{2}(t) & = & 2[(y^{2})^{2}(t)-(y^{1})^{2}(t)],\\[0.1cm]
D_{t}^{q}{x}^{3}(t) & = & - 2 (y^{2})^{2}(t),\\[0.1cm]
D_{t}^{q}{y}^{1}(t) & = & y^{1}(t)( x^{2}(t)-x^{1}(t)) \\[0.1cm]
D_{t}^{q}{y}^{2}(t) & = & y^{2}(t) ( x^{3}(t)- x^{2}(t))+b y^{2}(t), \label{(2.3)}
\end{array}\right.
\end{equation}
where $~a,b\in {\bf R^{\ast}} $ is a  parameter and $~q \in (0,1).$
\begin{Rem}{\rm  If we substitute $~q=1~$ and $~a=0~$ in the fractional model $~(2.3),~$ then the $~3-$dimensional Toda lattice is obtained. This dynamics was treated in paper \cite{puwa} from the point of view of Poisson geometry.} \hfill$\Box$
\end{Rem}
\markboth{M. Ivan}{On the 3-dimensional fractional-order Toda lattice with two controls}

Using the transformations $ x^{i} = x^{i}~$ for $~i=\overline{1,3},~ y^{1} = x^{4},~ y^{2}= x^{5},$ the system $(2.3)$ becomes:\\[-0.2cm]
\begin{equation}
\left\{ \begin{array} {lll}
D_{t}^{q}{x}^{1}(t) & = & 2 (x^{4})^{2}(t) + a x^{1}(t),\\[0.1cm]
D_{t}^{q}{x}^{2}(t) & = & 2[(x^{5})^{2}(t)-(x^{4})^{2}(t)],\\[0.1cm]
D_{t}^{q}{x}^{3}(t) & = & - 2 (x^{5})^{2}(t),~~~~~~~~~~~~~~~~~~~~~~ a,b \in {\bf R^{\ast}}, ~~~q \in (0,1),\\[0.1cm]
D_{t}^{q}{x}^{4}(t) & = & x^{4}(t)( x^{2}(t)-x^{1}(t)), \\[0.1cm]
D_{t}^{q}{x}^{5}(t) & = & x^{5}(t) ( x^{3}(t)- x^{2}(t))+b x^{5}(t). \label{(2.4)}
\end{array}\right.
\end{equation}

The fractional system $~(2.4)~$ is called the {\it 3-dimensional fractional-order Toda lattice with two controls around axes $~Ox^{1}~$ and $~Ox^{5}.~$ }

The initial value problem of the fractional system
$(2.4)$ can be represented in the following matrix form:\\[-0.2cm]
\begin{equation}
D_{t}^{\alpha}x(t)  =  x^{4}(t) A_{1} x(t) + x^{5}(t) A_{2} x(t) + A_{3} x(t) ,~~~~~~~~ x(0) =
x_{0},\label{(2.5)}
\end{equation}
where $0 < q < 1,~ x(t)= ( x^{1}(t),
 x^{2}(t), x^{3}(t), x^{4}(t), x^{5}(t))^{T}, ~t\in(0,\tau)$ and\\[-0.1cm]
\[
A_{1} = \left ( \begin{array}{ccccc}
0 & 0 & 0 & 2 & 0\\
0 & 0 & 0 & -2 & 0\\
0 & 0 & 0 & 0 & 0\\
-1 & 1 & 0 & 0 & 0\\
0 & 0 & 0 & 0 & 0\\
\end{array}\right ),~  A_{2} = \left ( \begin{array}{ccccc}
0 & 0 & 0 & 0 & 0\\
0 & 0 & 0 & 0 & 2\\
0 & 0 & 0 & 0 & -2\\
0 & 0 & 0 & 0 & 0\\
0 & -1 & 1 & 0 & 0\\
\end{array}\right ),~ A_{3} = \left ( \begin{array}{ccccc}
a & 0 & 0 & 0 & 0\\
0 & 0 & 0 & 0 & 0\\
0 & 0 & 0 & 0 & 0\\
0 & 0 & 0 & 0 & 0\\
0 & 0 & 0 & 0 & b\\
\end{array}\right ).
\]
\begin{Prop}
The initial value problem of the 3-dimensional fractional-order Toda lattice with two controls
$(2.5)$ has a unique solution.
\end{Prop}
{\it Proof.} Let $ f(x(t))= x^{4}(t) A_{1} x(t) + x^{5}(t) A_{2} x(t) + A_{3} x(t). $ It is
obviously continuous and bounded on $ D =\{ x \in {\bf R}^{5} |~
x^{i}\in [x_{0}^{i} - \delta, x_{0}^{i} + \delta]\}, i=\overline{1,5} $ for any
$\delta>0. $\\
 We have $~f(x(t)) - f(y(t)) = x^{4}(t) A_{1} x(t) + x^{5}(t) A_{2} x(t) + A_{3} x(t)-  y^{4}(t) A_{1} y(t) + y^{5}(t) A_{2} y(t) - A_{3} y(t) = g_{1}(t)+ g_{2}(t)+ g_{3}(t), $
where $~g_{1}(t)= x^{4}(t) A_{1} x(t) - y^{4}(t) A_{1} y(t),~ g_{2}(t)= x^{5}(t) A_{1} x(t) - y^{5}(t) A_{2} y(t)~$ and $~g_{3}(t)= A_{3} x(t) - A_{3} y(t).~$ Then\\[0.1cm]
 $(a)~~|f(x(t)) - f(y(t))|\leq |g_{1}(t)| + |g_{2}(t)| + |g_{3}(t)|. $\\
Using reasoning analogous to that in the proof of the Proposition 2.1 in \cite{ivmi}, we can show that:\\[0.1cm]
$(b)~~~|g_{1}(t)| \leq (\|A_{1}\|+ |y^{4}(t)|)\cdot |x(t)- y(t)|,~~~|g_{2}(t)| \leq (\|A_{2}\|+ |y^{5}(t)|)\cdot |x(t)- y(t)|,$   and  $~|g_{3}(t)| \leq \|A_{3}\| \cdot |x(t)- y(t)|.$\\[0.1cm]
According to $(b)~$ the relation $(a)$ becomes\\[0.1cm]
$(c)~~~|f(x(t)) - f(y(t))|\leq  (\|A_{1}\| + \|A_{2}\| + \|A_{3}\| +|y^{4}(t)| + |y^{5}(t)|)\cdot |x(t)-y(t)|.$\\[0.1cm]
Replacing  $~\|A_{1}\|= \|A_{2}\|=\sqrt{10},~ \|A_{3}\| =2 |a|~$ and using the inequalities $~|y^{i}(t)|\leq |x_{0}| + \delta, ~i=1,5~$ from the
relation $ (c), $ we deduce that\\[0.1cm]
$(d)~~~|f(x(t)) - f(y(t))|\leq  L\cdot |x(t)-y(t)|,~~~~~
\hbox{where}~ L =2\sqrt{10} + \sqrt{a^{2}+b^{2}} + 3(|x_{0}| + \delta) > 0.$\\[0.1cm]
 The inequality $(d)$ shows that $ f(x(t))$
satisfies a Lipschitz condition. Based on the results of Theorems
$1$ and $2$ in \cite{diet}, we can conclude that the initial value
problem of the system $(2.5)$ has a unique solution. \hfill$\Box$

For the fractional system $ (2.4) $ we introduce the following notations:\\[-0.2cm]
 \begin{equation}
 \left \{\begin{array}{l}
f_{1}(x)=2 (x^{4})^{2} + a x^{1},~~~ f_{2}(x)=2[(x^{5})^{2}-(x^{4})^{2}],~~~ f_{3}(x)=- 2 (x^{5})^{2},\\
f_{4}(x)=x^{4}( x^{2}-x^{1}),~~~~ f_{5}(x)=x^{5}( x^{3}- x^{2})+b x^{5}.\label{(2.6)}
\end{array}\right.
\end{equation}
\begin{Prop}
{\it The equilibrium states of the 3-dimensional fractional-order Toda lattice
$(2.4)$ are given as the following family}:\\
$ E:=\{e_{0}= (0,0,0,0,0)\} \cup E_{1},~$ where\\
$~ E_{1}=\{ e_{23}^{k,m} =(0, k, m, 0, 0)\in {\bf R}^{5} |~ k,m \in {\bf R},~ k\neq 0~ or~ m\neq 0\}.$
\end{Prop}
{\it Proof.} The equilibrium states are solutions of the equations
$~f_{i}(x)=0, i=\overline{1,5}$ where $~f_{i},~i=\overline{1,5}$
are given by (2.6).\hfill$\Box$

The set $~E_{1}~$ of the equilibrium states of the system $~(2.3)~$ includes the following three remarkable subsets:\\
$E_{1}^{1}=\{ e_{2}^{k}:=e_{23}^{k,0} =(0, k, 0, 0, 0)\in {\bf R}^{5} |~ k \in {\bf R}^{\ast}\},~$\\
$E_{1}^{2}=\{ e_{3}^{m}:=e_{23}^{0,m} =(0, 0, m, 0, 0)\in {\bf R}^{5} |~ m \in {\bf R}^{\ast}\}~$ and\\
$E_{1}^{3}=\{ e_{23}^{m,m} =(0, m, m, 0, 0)\in {\bf R}^{5} |~ m \in {\bf R}^{\ast}\}.$

Let us we present the study of asymptotic stability of equilibrium states for the fractional system $(2.4). $ For this study we apply
the Matignon's test.

The Jacobian matrix associated to system $(2.4)$ is:
\[
J(x,a,b) = \left ( \begin{array}{ccccc}
a & 0 & 0 & 4x^{4} & 0 \\
0 & 0 & 0 & -4x^{4} & 4x^{5} \\
0 & 0 & 0 & 0 & -4x^{5} \\
-x^{4} & x^{4} & 0 & x^{2}-x^{1} & 0 \\
0 & -x^{5} & x^{5} & 0 & x^{3}-x^{2}+b \\
\end{array}\right ).\\[-0.1cm]
\]
\begin{Prop} {\rm (\cite{mati})}
Let $ x_{e} $ be an equilibrium state of fractional differential system $(2.4)$ and
$ J(x_{e},a,b) $ be the Jacobian matrix $J(x,a,b)$ evaluated at $ x_{e}$.

 $(i)~ x_{e}$ is locally asymptotically stable, if and only if  all eigenvalues
$ \lambda(J(x_{e},a,b)) $ of  $ J(x_{e},a,b) $ satisfy:
\begin{equation}
| arg(\lambda (J(x_{e},a,b))) | > \displaystyle\frac{q\pi}{2}.\label{(2.7)}
\end{equation}
 $(ii)~ x_{e} $ is locally stable, if and only if either it is asymptotically stable, or the
critical eigenvalues satisfying $~| arg(\lambda (J(x_{e},a,b)) | = \displaystyle\frac{q \pi}{2}~$ have geometric
multiplicity one.\hfill$\Box$
\end{Prop}
Using the notation: $~~\tilde{q}:=\frac{2}{\pi}|arg(\lambda (J(x_{e},a,b)))|~$ and applying Proposition $ 2.3~$ one obtains the following corollaries.
\begin{Cor}{\rm (\cite{gi22})}
$(i)~$ The equilibrium state $~x_{e}~$  of the fractional model $ (2.4) $ is asymptotically stable if and only if the difference $~q-\tilde{q}~$ is strictly negative. More precisely, $~x_{e}~$ is asymptotically stable for all $~q\in (0, \tilde{q}).~$

$(ii)~$ If $~q-\tilde{q} > 0,~$ then  $~x_{e}~$ is unstable and the fractional model $ (2.4) $  may exhibit chaotic behavior. More precisely, $~x_{e}~$ is unstable $~(\forall) q\in (\tilde{q}, 1).~$ \hfill$\Box$
\end{Cor}
\begin{Cor} {\rm (\cite{gi22})} Let $~x_{e}~$ be an equilibrium state of the fractional model $~(2.4)~$ and $~\lambda_{i},~i=\overline{1,5}~$ the eigenvalues of $~J(x_{e},a,b).$

$(i)~$ If one of the eigenvalues $~\lambda_{i},~i=\overline{1,5}~$ is equal to zero or it is positive, then $~x_{e}~$ is unstable  for all $~q\in (0,1).$

$(ii)~$ If $~\lambda_{i} < 0, $ for all $~i=\overline{1,5},~$ then  $~x_{e}~$   is asymptotically stable  $~(\forall)~q\in (0,1).$
\end{Cor}
\begin{Prop}
The equilibrium states $~e_{0}~$ and $~ e_{23}^{k,m}\in E_{1}~$ are unstable $ (\forall) q \in (0,1).$
\end{Prop}
{\it Proof.} The characteristic polynomial of  $~
J(e_{23}^{k,m},a,b) =\left (\begin{array}{ccccc}
  a & 0 & 0 & 0 & 0 \\
0 & 0 & 0 & 0 & 0 \\
0 & 0 & 0 & 0 & 0 \\
0 & 0 & 0 & k & 0 \\
0 & 0 & 0 & 0 & m-k+b\\
\end{array}\right )~$
is\\
 $~ p_{J(e_{23}^{k,m},a,b)}(\lambda) = \det ( J(e_{23}^{k,m},a,b) -
\lambda I) = - \lambda^{2} (\lambda -a)(\lambda -k)(\lambda -b +k-m).~$
 The equation $ ~p_{J(e_{23}^{k,m},a,b)}(\lambda) = 0 $  has the root $ \lambda_{1} = 0.~$ By Corollary  2.2(i), follows that $~e_{0}~$ and $~ e_{23}^{k,m}\in E_{1}~$ are unstable for all $ q\in (0,1).$ \hfill$\Box$
\begin{Rem}{\rm  Toda-type dynamical systems have been studied from various research directions by many authors. More specifically, from the point of view of Poisson geometry, Toda lattices were discussed in the papers \cite{dami, pudo, puwa}. Fractional 
dynamical systems associated to Toda-type dynamics and to Volterra lattices have been investigated in  \cite{gi22,igom, miva23, imgo}.} \hfill$\Box$
\end{Rem}

\section{Controllability of chaotic behaviors of the $ 3-$dimensional fractional-order Toda lattice with two controls $~(2.4) $}

In this section we will discuss how to stabilize the unstable equilibrium states of the fractional system $ (2.4) $ via fractional-order derivative. For this purpose we will apply the general  method for to control the stability  of $~(2.4)~$ at an equilibrium point \cite{gi22}.

Let  $ x_{e}$ be an unstable equilibrium point of the $ 3-$dimensional fractional-order Toda lattice  with two controls $~(3.2).~$ We associate to $(2.4)$ a new fractional-order system with (external) controls and given by:\\[-0.2cm]
\begin{equation}
\left\{ \begin{array} {lll}
D_{t}^{q}{x}^{1}(t) & = & 2 (x^{4})^{2}(t) + a x^{1}(t)+ u_{1}(t),\\[0.1cm]
D_{t}^{q}{x}^{2}(t) & = & 2[(x^{5})^{2}(t)-(x^{4})^{2}(t)] + u_{2}(t),\\[0.1cm]
D_{t}^{q}{x}^{3}(t) & = & - 2 (x^{5})^{2}(t)+ u_{3}(t),~~~~~~~~~~~~~~~~~~~~~~ a,b \in {\bf R^{\ast}}, ~~~q \in (0,1),\\[0.1cm]
D_{t}^{q}{x}^{4}(t) & = & x^{4}(t)( x^{2}(t)-x^{1}(t))+ u_{4}(t), \\[0.1cm]
D_{t}^{q}{x}^{5}(t) & = & x^{5}(t) ( x^{3}(t)- x^{2}(t))+b x^{5}(t) + u_{5}(t), \label{(3.1)}
\end{array}\right.
\end{equation}
where $ u_{i}(t), i=\overline{1,5}$ are control functions.

We take the control functions $ u_{i}(t), i=\overline{1,5}, $ given by:\\[-0.2cm]
 \begin{equation}
u_{1}= 0,~~u_{2} = c_{1} (x^{2}- x_{e}^{2}),~~ u_{3} = c_{2} (x^{3}- x_{e}^{3}),~~ u_{4}  = c_{3} (x^{4}- x_{e}^{4}),~~ u_{5}= 0,\label{(3.2)}
\end{equation}
where $~ c_{1}, c_{2}, c_{3} \in {\bf R}^{*}.$

With the control functions $(3.2), $ the system $(3.1) $ becomes:\\[-0.2cm]
\begin{equation}
\left\{ \begin{array} {lll}
D_{t}^{q}{x}^{1}(t) & = & 2 (x^{4})^{2}(t) + a x^{1}(t),\\[0.1cm]
D_{t}^{q}{x}^{2}(t) & = & 2[(x^{5})^{2}(t)-(x^{4})^{2}(t)] + c_{1} (x^{2}(t)- x_{e}^{2}),\\[0.1cm]
D_{t}^{q}{x}^{3}(t) & = & - 2 (x^{5})^{2}(t)+ c_{2} (x^{3}(t)- x_{e}^{3}),\\[0.1cm]
D_{t}^{q}{x}^{4}(t) & = & x^{4}(t)( x^{2}(t)-x^{1}(t))+c_{3} (x^{4}(t)- x_{e}^{4}), \\[0.1cm]
D_{t}^{q}{x}^{5}(t) & = & x^{5}(t) ( x^{3}(t)- x^{2}(t))+b x^{5}(t), \label{(3.3)}
\end{array}\right.
\end{equation}
where  $~a, b, c_{1}, c_{2}, c_{3} \in{\bf R}^{\ast}~$ are control parameters and $~q \in (0,1).$

The fractional system $ (3.3) $ is called the {\it controlled fractional-order Toda lattice associated  to
$~(3.1) $ at  $~x_{e}.$}

If one selects the appropriate parameters $~a,b, c_{1}, c_{2}, c_{3}\in {\bf R}^{\ast} $  which then make the eigenvalues of the linearized equation of $(3.3)$ satisfy one of the conditions from Proposition 2.3, then the trajectories of $ (3.3) $ asymptotically approaches the unstable equilibrium state $x_{e}$ in the sense that $\lim_{t\rightarrow \infty} \|x(t)-x_{e}\|= 0$, where $\|\cdot\|$ is the Euclidean norm. In the case when the equilibrium state $~x_{e}~$ is unstable, then fractional model $ (3.3) $  may exhibit chaotic behavior.

The Jacobian matrix of the controlled fractional model $(3.3)$  is:\\[-0.2cm]
\[
J(x,a, b, c_{1}, c_{2}, c_{3}) = \left ( \begin{array}{ccccc}
a & 0 & 0 & 4x^{4} & 0 \\
0 & c_{1} & 0 & -4x^{4} & 4x^{5} \\
0 & 0 & c_{2} & 0 & -4x^{5} \\
-x^{4} & x^{4} & 0 & x^{2}-x^{1}+c_{3} & 0 \\
0 & -x^{5} & x^{5} & 0 & x^{3}-x^{2}+b \\
\end{array}\right ).\\[-0.1cm]
\]
\begin{Prop}
Let be the controlled fractional system  $(3.3)~$ and $~e_{0}= (0,0,0,0,0).$\\
$~~~~~~~$ {\bf 1.} $~a < 0.$\\
$(i)~$ If $~ b<0, c_{1} <0,~c_{2} <0~$  and $~c_{3} <0,~$ then $ e_{0}~$ is asymptotically stable $(\forall) q\in (0,1).$\\
$(ii)~$  If at least one of the parameters $~b, c_{1}, c_{2} ~$ and $~c_{3} $ is positive, then $ e_{0}~$ is unstable $(\forall) q\in (0,1).$\\
$~~~~~~~$ {\bf 2.} $~a > 0.~$ If  $~b, c_{1}, c_{2}, c_{3}\in {\bf R}^{\ast},~$ then $~ e_{0} $ is unstable  $ (\forall) q\in (0,1).$
\end{Prop}
{\it Proof.} The characteristic polynomial of the Jacobian matrix $~J_{0}:=J(e_{0},a, b, c_{1}, c_{2}, c_{3})~$
 is $~ p_{0}(\lambda) = \det ( J_{0} - \lambda I) = (\lambda -a)(\lambda - b) (\lambda -c_{1})(\lambda -c_{2})(\lambda -c_{3}).~$
The roots of equation $~ p_{0}(\lambda) = 0~$ are $\lambda_{1}= a,~\lambda_{2}= b,~\lambda_{j+2}=c_{j},~j=\overline{1,3}.$\\
{\bf 1.} {\bf Case} $~a <0~$ and $~q\in (0,1).~$  Then $~\lambda_{1} <0.~$\\
$(i)~$ We have $~\lambda_{j+1}<0~$ for $~j=\overline{1,4}~$ if and only if $~b<0, c_{1}<0, c_{2}<0~$  and  $~c_{3}<0.~$ Then $~\lambda_{i} <0, i=\overline{1,5}~$ and according to Corollary 2.2(ii), $~ e_{0}~$ is asymptotically stable.\\
$(ii)~$ We suppose $~ a < 0~$  and $~b>0,~$ or $~c_{1}>0~$ or $~c_{2}>0~$ or $~c_{3}>0.~$ Then $~J_{0}~$ has at least a positive eigenvalue and by Corollary  2.2(i), $~ e_{0} $ is unstable.\\
{\bf 2.} {\bf Case} $~a >0~$ and $~q\in (0,1).~$ Since $~\lambda_{1}>0,~ J_{0}~$ has at least a positive eigenvalue $~(\forall) b, c_{1}, c_{2}, c_{3}\in {\bf R}^{\ast}.~$ Then  $ e_{0} $ is unstable.\hfill$\Box$
\begin{Prop}
Let  be the controlled fractional system  $(3.3).~$ The equilibrium states $ e_{23}^{k,m}=(0, k, m, 0, 0)\in E_{1}~$  have the following behavior:\\
$~~~~~~~$ {\bf 1.} $~a < 0~$ and $~q\in (0,1).$\\
$(i)~$ If $~ c_{1} < 0~$ and $~c_{2} < 0,~$ then $ e_{23}^{k,m}~$ is asymptotically stable for all $~k\in (-\infty, -c_{3})~$ and $~m\in (-\infty, k-b).~$\\
$(ii)~$ If $~c_{1} < 0~$ and $~c_{2} <0,~$ then $~ e_{23}^{k,m}~$ is unstable for all $~k\in (-c_{3}, \infty)~$ or $~m\in (k-b, \infty).$\\
$(iii)~$ If $~c_{1} > 0~$ or $~c_{2} > 0~$ and $~b, c_{3} \in {\bf R}^{\ast},~$ then $~ e_{23}^{k,m}~$ is unstable for all $~k \in {\bf R}^{\ast}~$ or $~m \in {\bf R}^{\ast}.$\\
$~~~~~~~$ {\bf 2.} $~a> 0~$ and $~q\in (0,1).~$ If $~b,c_{1}, c_{2}, c_{3}\in {\bf R}^{\ast},~$ then $~ e_{23}^{k,m}~$ is unstable for all $~k \in {\bf R}^{\ast}~$ or $~m \in {\bf R}^{\ast}.$
\end{Prop}
{\it Proof.} Denote $~J_{1}:= J(e_{23}^{k,m}, a, b, c_{1}, c_{2}, c_{3}).~$ The characteristic polynomial  of the matrix
\begin{center}
$~J_{1} = \left ( \begin{array}{ccccc}
a & 0 & 0 & 0 & 0 \\
0 & c_{1} & 0 & 0 & 0 \\
0 & 0 & c_{2} & 0 & 0 \\
0 & 0 & 0 & k + c_{3} & 0 \\
0 & 0 & 0 & 0 & m- k+ b \\
\end{array}\right )~$ is
\end{center}
$~p_{1}(\lambda) = \det ( J_{1} -\lambda I) = (\lambda -a)(\lambda -c_{1})(\lambda -c_{2}) (\lambda - c_{3} - k)(\lambda -b+ k- m).~$
The roots of equation $~ p_{1}(\lambda) = 0~$ are $~\lambda_{1}= a,~ \lambda_{2}= c_{1},~\lambda_{3}= c_{2},~\lambda_{4}=c_{3}+k,~\lambda_{5}=b-k+m.$\\
{\bf 1.} {\bf Case} $~a <0~$ and $~q\in (0,1).~$  Then $~\lambda_{1} <0.~$\\
$(i)~$ We suppose $~c_{1}<0~$ and $~c_{2}<0. $ Then $~\lambda_{2}<0~$ and $~\lambda_{3}<0.~$ We have $~\lambda_{i} <0, i=\overline{1,5}~$ if and only if $~k\in (-\infty, -c_{3})~$ and $~m\in (-\infty, k-b).~$ By Corollary 2.2(ii),$~ e_{23}^{k,m}~$ is asymptotically stable.\\
$(ii)~$ We suppose $~c_{1}<0~$ and $~c_{2}<0.~$ For $~k\in (-c_{3}, \infty)~$ it follows that $~\lambda_{4} >0.~$ Also, for  $~m\in (k-b, \infty),~$ it follows that $~\lambda_{5} >0.~$ In the both  cases, the matrix  $~J_{1}~$ has  a positive eigenvalue. According to Corollary  2.2(i),$~ e_{23}^{k,m}~$ is unstable for all $~k\in (-c_{3}, \infty)~$ or $~m\in (k-b, \infty).$\\
$(iii)~$ We suppose $~c_{1}>0~$ or $~ c_{2}>0~$ and $~b, c_{3} \in {\bf R}^{\ast}.~$  Then $~\lambda_{2}>0~$ or $~\lambda_{3}>0~$ or and $~J_{1}~$ has at least a positive eigenvalue. By Corollary  2.2(i),$~ e_{23}^{k,m}~$ is unstable  $~(\forall) k \in {\bf R}^{\ast}~$ or $~m \in {\bf R}^{\ast}.$\\
{\bf 2.} {\bf Case} $~a >0~$ and $~q\in (0,1).~$  Since $~\lambda_{1} >0,~$ the matrix $~J_{1}~$ has  a positive eigenvalue $~(\forall) b, c_{1}, c_{2}, c_{3} \in {\bf R}^{\ast}.~$ Then $~ e_{23}^{k,m}~$ is unstable $~(\forall) k \in {\bf R}^{\ast}~$ or $~m \in {\bf R}^{\ast}.$\hfill$\Box$

Using Proposition $~3.2,~$ it is easy to prove the following three corollaries.
\begin{Cor}
Let  be the controlled fractional system  $(3.3).~$ The equilibrium states $ e_{2}^{k}=(0, k, 0, 0, 0)\in E_{1}^{1}~$  have the following behavior:\\
$~~~~~~~$ {\bf 1.} $~a < 0~$ and $~q\in (0,1).$\\
$(i)~$ If $~ c_{1} < 0, c_{2} < 0, b <-c_{3},$ then $ e_{2}^{k}~$ is asymptotically stable $~(\forall) k\in (b, -c_{3}).$\\
$(ii)~$ If $~c_{1} < 0, c_{2} < 0, b >-c_{3}, $ then $ e_{2}^{k}~$ is unstable $~(\forall) k\in (-c_{3}, \infty)~$ or $~k\in (-\infty, b).$\\
$(iii)~$ If $~c_{1} > 0~$ or $~c_{2} > 0~$ and $~b, c_{3} \in {\bf R}^{\ast},~$ then $~ e_{2}^{k}~$ is unstable $~(\forall) k \in {\bf R}^{\ast}.$\\
$~~~~~~~$ {\bf 2.} $~a> 0~$ and $~q\in (0,1).~$ If $~b,c_{1}, c_{2}, c_{3}\in {\bf R}^{\ast},~$ then $~ e_{2}^{k}~$ is unstable $~(\forall) k \in {\bf R}^{\ast}.$
\end{Cor}
\begin{Cor}
Let  be the controlled fractional system  $(3.3).~$ The equilibrium states $ e_{3}^{m}=(0, 0, m, 0, 0)\in E_{1}^{2}~$  have the following behavior:\\
$~~~~~~~$ {\bf 1.} $~a < 0~$ and $~q\in (0,1).$\\
$(i)~$ If $~ c_{1} < 0, c_{2} < 0, c_{3}<0,$ then $ e_{3}^{m}~$ is asymptotically stable $~(\forall) m\in (-\infty, -b).$\\
$(ii)~$ If $~c_{1} < 0, c_{2} < 0, $ then $ e_{3}^{m}~$ is unstable for $~c_{3}>0~$ or $~m\in (-b, \infty).$\\
$(iii)~$ If $~c_{1} > 0~$ or $~c_{2} > 0~$ and $~b, c_{3} \in {\bf R}^{\ast},~$ then $~ e_{3}^{m}~$ is unstable $~(\forall) m \in {\bf R}^{\ast}.$\\
$~~~~~~~$ {\bf 2.} $~a> 0~$ and $~q\in (0,1).~$ If $~b,c_{1}, c_{2}, c_{3}\in {\bf R}^{\ast},~$ then $~ e_{3}^{m}~$ is unstable $~(\forall) m \in {\bf R}^{\ast}.$
\end{Cor}
\begin{Cor}
Let  be the controlled fractional system  $(3.3).~$ The equilibrium states $ e_{23}^{m,m}=(0, m, m, 0, 0)\in E_{1}^{3}~$  have the following behavior:\\
$~~~~~~~$ {\bf 1.} $~a < 0~$ and $~q\in (0,1).$\\
$(i)~$ If $~ c_{1} < 0, c_{2} < 0, b<0,$ then $ e_{23}^{m,m}~$ is asymptotically stable $~(\forall) m\in (-\infty, -c_{3}).$\\
$(ii)~$ If $~c_{1} < 0, c_{2} < 0, $ then $ e_{23}^{m,m}~$ is unstable for $~b>0~$ or $~m\in (-c_{3}, \infty).$\\
$(iii)~$ If $~c_{1} > 0~$ or $~c_{2} > 0~$ and $~b, c_{3} \in {\bf R}^{\ast},~$ then $~ e_{23}^{m,m}~$ is unstable $~(\forall) m \in {\bf R}^{\ast}.$\\
$~~~~~~~$ {\bf 2.} $~a> 0~$ and $~q\in (0,1).~$ If $~b,c_{1}, c_{2}, c_{3}\in {\bf R}^{\ast},~$ then $~ e_{23}^{m,m}~$ is unstable $~(\forall) m \in {\bf R}^{\ast}.$
\end{Cor}
\begin{Ex}
{\rm  Let be the controlled fractional Toda lattice $ (3.3)~$ and $~q\in (0, 1).$\\
$(1)~$ We select $ a=-0.8, b=-0.2, c_{1}=-0.03, c_{2}=-0.02~$ and  $~c_{3}=-0.001.~$ By Proposition  $~3.1.1(i),~e_{0}=(0,0,0,0,0)~$ is asymptotically stable.\\
$(2)~$ Let $ a=-0.8, b=0.2, c_{1}=-0.03, c_{2}=-0.02~$ and  $~c_{3}=-0.001.~$ According to Proposition  $~3.1.1(ii),~e_{0}~$ is unstable. In other words, in this case the fractional model $~(3.3)~$ behaves chaotically around the equilibrium point $~e_{0}.$}\\[-0.4cm]
\end{Ex}
\begin{Ex}
{\rm  Let be the controlled fractional Toda lattice $ (3.3)~$ and $~q\in (0, 1).$\\
$(1)~$ Let $ a =-0.45, b=1, c_{1}=-0.2, c_{2}=-0.15, c_{3}=1.01,~k=1.~$ By Proposition  $~3.2.1(i),~ e_{23}^{1,m}=(0, 1, m,0,0)~$ is asymptotically stable $~(\forall) m\in (0, \infty).$ \\
$(2)~$ Let $ a =-0.45, b =1, c_{1}=-0.2, c_{2}=-0.15,~c_{3}=1.01,~k=-1.~$ Applying Proposition $~3.3.1(ii),~e_{23}^{-1,m}=(0,-1,m,0,0)~$ is unstable $~(\forall) m\in {\bf R}^{\ast}.$}
\end{Ex}
\begin{Ex}
{\rm  Let be the controlled fractional Toda lattice $ (3.3)~$ and $~q\in (0, 1).$\\
$(1)~$ Let $ a =-1, b=0.32, c_{1}=-0.25, c_{2}=-0.12, c_{3}=-1.05.~$ By Corollary $~3.1.1(i),~ e_{2}^{k}=(0, k, 0,0,0)~$ is asymptotically stable $~(\forall) k\in (0.32, 1.05).$\\
$(2)~$ Let $~a =-1, b=1.2, c_{1}=-0.25, c_{2}=-0.12, c_{3}=-1.05.~$ Applying Corollary $~3.1.1(ii),~e_{2}^{k}=(0,k,0,0,0)~$ is unstable for all $~k\in {\bf R}^{\ast}.$}
\end{Ex}
\begin{Ex}
{\rm  Let be the controlled fractional Toda lattice $ (3.3)~$ and $~q\in (0, 1).$\\
$(1)~$ Let $ a =-0.9, b=0.08, c_{1}=-0.4, c_{2}=-0.22, c_{3}=-0.06.~$ By Corollary $~3.2.1(i),~e_{3}^{m}=(0, 0, m,0,0)~$ is asymptotically stable $~(\forall) m\in (-\infty, -0.08).$\\
$(2)~$ Let $ a =-0.9, b=1, c_{1}=-0.4, c_{2}=-0.22, c_{3}=-0.06.~$ Applying Corollary $~3.2.1(ii),~e_{3}^{m}=(0,0,m,0,0)~$ is unstable for all $~m\in (-1, \infty)\setminus \{0\}.$}
\end{Ex}
\begin{Ex}
{\rm  Let be the controlled fractional Toda lattice $ (3.3)~$ and $~q\in (0, 1).$\\
$(1)~$ Let $ a =-0.75, b=-0.81, c_{1}=-0.36, c_{2}=-0.01, c_{3}=0.48.~$ By Corollary $~3.3.1(i),~e_{23}^{m,m}=(0, m, m,0,0)~$ is asymptotically stable $~(\forall) m\in (-\infty, -0.48).$\\
$(2)~$ Let $ a =-0.75, b=-0.81, c_{1}=-0.36, c_{2}=-0.01, c_{3}=0.48.~$ Applying Corollary $~3.3.1(ii),~e_{23}^{m,m}=(0,m,m,0,0)~$ is unstable for all $~m\in (-0,48, \infty)\setminus \{0\}.$}\\[-0.5cm]
\end{Ex}
\begin{Rem}{\rm The general method for controlling the stability of a fractional dynamical system has also been used for to study the stability of other fractional systems, for example in 
 \cite{mi22, mihi, miva23, miva24, mivan24}.} \hfill$\Box$
\end{Rem}

\section{Numerical integration of the controlled fractional-order Toda lattice $ (3.3)$}

In this section we apply the fractional Euler's method (FEM) to numerically integrate the controlled fractional-order Toda lattice $ (3.3).~$ For the description and application of FEM's can be consulted \cite{aheg, mi22, odmo}.

Consider the following general form of the initial value problem (IVP) with Caputo derivative \cite{odmo}:\\[-0.4cm]
\begin{equation}
D_{t}^{q} y(t) = f(t,y(t)),~~~ y(0)=y_{0},~~~t\in I=[0,T],~T>0  \label{(4.1)}
\end{equation}
where $~y: I \rightarrow {\bf R}^{n},~f: {\bf R}^{n} \rightarrow {\bf R}^{n}~$ is a continuous nonlinear function and
$ q\in (0,1).$

Every solution of the initial value problem given by $ (4.1) $ is also a solution of the following {\it Volterra fractional integral equation}:\\[-0.4cm]
\begin{equation}
y(t) = y(0) +~ I_{t}^{q}f(t,y(t)), \label{(4.2)}
\end{equation}
where $ I_{t}^{q} $ is the $ q-$order Riemann-Liouville integral operator. Moreover, every solution of $ (4.2) $ is a solution of the (IVP) $ (4.1).$

To integrate the fractional equation $( 4.1),$ means to find the solution of $ (4.2) $ over the interval $~[0,T]. $ In this context, a set of points $~(t_{j}, y(t_{j})) $ are produced which are used as approximated values. The interval $ [0,T] $  is partitioned into $ n $ subintervals $~[t_{j}, t_{j+1} ] $  each equal width $~ h =\frac{T}{n}, ~ t_{j} = j h $  for  $~ j = 0,1,..., n.~$ It computes an approximation denoted as $~y_{j+1} ~$ for $~ y ( t_{j+1}),~ j= 0, 1, ... .$

The {\it general formula of the fractional Euler's method} for to compute the elements $~y_{j},$ is\\[-0.4cm]
\begin{equation}
y_{j+1} = y_{j} + \displaystyle\frac{h^{q}}{\Gamma(q+1)} f( t_{j}, y(t_{j})),~~~~~ t_{j+1} = t_{j} + h,~~~ j=0, 1, ..., n. \label{(4.3)}
\end{equation}
For more details, see \cite{ ahme, danc, odmo}.

We will now apply the above considerations to the controlled fractional-order Toda lattice $ (3.3).~$ For this, consider the following fractional differential equations\\[-0.2cm]
\begin{equation}
\left\{\begin{array}{lcl}
 D_{t}^{q} x^{i}(t) & = &
F_{i}(x^{1}(t), x^{2}(t), x^{3}(t), x^{4}(t), x^{5}(t)),~~~ i=\overline{1,5}, ~~t\in
(t_{0}, \tau),~q \in (0,1)\\
x(t_{0}) &=& (x^{1}(t_{0}), x^{2}(t_{0}), x^{3}(t_{0}), , x^{4}(t_{0}) , x^{5}(t_{0}))
\end{array}\right.\label{(4.4)}\\[-0.2cm]
\end{equation}
where\\[-0.4cm]
\begin{equation}
\left\{\begin{array}{lcl} F_{1}(x(t)) & = & 2 (x^{4})^{2}(t) + a x^{1}(t),\\
F_{2}(x(t)) & = & 2[(x^{5})^{2}(t)-(x^{4})^{2}(t)] + c_{1} (x^{2}(t)- x_{e}^{2}),\\
F_{3}(x(t)) & = & - 2 (x^{5})^{2}(t)+ c_{2} (x^{3}(t)- x_{e}^{3}),\\
F_{4}(x(t)) & = & x^{4}(t)( x^{2}(t)-x^{1}(t))+c_{3} (x^{4}(t)- x_{e}^{4}),\\
F_{5}(x(t)) & = & x^{5}(t) ( x^{3}(t)- x^{2}(t))+b x^{5}(t),\\
\end{array}\right.\label{(4.5)}
\end{equation}
 where  $~a, b, c_{1}, c_{2}, c_{3}\in {\bf R}^{*}.$

 Since the functions $ F_{i}(x(t)), i=\overline{1,5} $ are continuous, the initial value problem $(4.4)$ is equivalent
 to system of Volterra integral equations, which is given as follows:\\[-0.2cm]
\begin{equation}
x^{i}(t)~=~ x^{i}(0)  + ~I_{t}^{q} F_{i}(x^{1}(t), x^{2}(t), x^{3}(t)), x^{4}(t)),x^{5}(t)),
~~~~~i=\overline{1,5}.\label{(4.6)}
\end{equation}

The system $(4.6)$ is called the {\it Volterra integral equations associated to controlled fractional-order Toda lattice} $(3.3)$.

The problem for solving the system $(4.4)$ is reduced to one of solving a sequence of systems of fractional equations in
increasing dimension on successive intervals $[j, (j+1)]$.

For the numerical integration of the system $(4.4)$ one can use the fractional Euler's method (the formula $(4.3)$ ), which is expressed as follows:\\[-0.4cm]
\begin{equation}
x^{i}(j+1)=x^{i}(j)+ \frac{h^{q}}{\Gamma (q+1)} F_{i}(x^{1}(j), x^{2}(j),
x^{3}(j)), x^{4}(j)), x^{5}(j)),~~~ i=\overline{1,5},\label{(4.7)}
\end{equation}
where $ j=0,1,2,...,N,  h=\displaystyle\frac{T}{N}, T>0, N>0.~$

More precisely, the numerical integration of the fractional system $(4.4)$ is given by:\\[-0.4cm]
\begin{equation}
\left \{ \begin{array}{ll} x^{1}(j+1) &= x^{1}(j)+
h^{q}~\displaystyle\frac{1}{\Gamma(q+1)} (2 (x^{4})^{2}(j) + a x^{1}(j))\\[0.3cm]
x^{2}(j+1) &= x^{2}(j)+
h^{q}~\displaystyle\frac{1}{\Gamma(q+1)}(2[(x^{5})^{2}(j)-(x^{4})^{2}(j)] + c_{1} (x^{2}(j)- x_{e}^{2}) )\\[0.3cm]
 x^{3}(j+1) &= x^{3}(j) +
h^{q}~\displaystyle\frac{1}{\Gamma(q+1)}(- 2 (x^{5})^{2}(j)+ c_{2} (x^{3}(j)- x_{e}^{3}) )\\[0.3cm]
x^{4}(j+1) &= x^{4}(j) +
h^{q}~\displaystyle\frac{1}{\Gamma(q+1)}(x^{4}(j)( x^{2}(j)-x^{1}(j))+c_{3} (x^{4}(j)- x_{e}^{4}) )\\[0.3cm]
x^{5}(j+1) &= x^{5}(j) +
h^{q}~\displaystyle\frac{1}{\Gamma(q+1)}(x^{5}(j) ( x^{3}(j)- x^{2}(j))+b x^{5}(j) )\\[0.3cm]
x^{i}(0)&= x_{e}^{i}+\varepsilon,~~~i=\overline{1,5}.\\
\end{array}\right. \label{(4.8)}
\end{equation}

Using \cite{diet, odmo}, we have that the numerical algorithm given by $(4.8)$ is convergent.\\[-0.4cm]
\begin{Ex}
{\rm  Let us we present the numerical integration of the controlled fractional-order Toda lattice with two controls  which has considered in Example $~3.2.1(i).$ For this we apply the algorithm $ (4.8) $ and software Maple. Then, in $(4.8)$ we take: $~ a =-0.45, b = 1, c_{1}=-0.2, c_{2}=-0.15, c_{3}=1.01,  k=1~$  and  $~m = 0.6. $ It is known that the equilibrium state $~ e_{23} = (0, 1, 0.6, 0,0)~$ is asymptotically stable for $~ q=0.8.$

For the numerical simulation of solution of the above fractional model we use the rutine {\it Maple. fract-Toda lattice-with-two controls}, denoted by [fr.Toda lattice-two controls]. Applying this program for $~h = 0.01, \varepsilon= 0.01, x^{1}(0)= \varepsilon, x^{2}(0)=1 + \varepsilon, x^{3}(0)= 0.6 + \varepsilon, x^{4}(0)= \varepsilon, x^{5}(0)= \varepsilon,~N = 100, t = 102 $ one obtain the orbits $~(n, x^{1}(n)), (n, x^{2}(n)), (n, x^{3}(n)), (n, x^{4}(n))~$ and $~(n, x^{5}(n)), $  for $ q = 0.8.$}\hfill$\Box$\\[-0.4cm]
\end{Ex}
Finally, we present the rutine  [fr.Toda lattice -two controls]:\\[0.2cm]
$ \# $ Fractional equations associated to controlled Toda lattice for q=0.8\\[0.1cm]
Dx1/dt=2*x4$\wedge$ 2 + a*x1; Dx2/dt=2*(x5$\wedge$ 2 - x4$\wedge$ 2)+c1*(x2-x2e); Dx3/dt=-2*x5$\wedge$ 2 +c2 *(x3-x3e); Dx4/dt=x4*(x2 -x1) +c3 *(x4-x4e); Dx5/dt=x5*(x3 -x2)+ b*x5;\\[0.2cm]
{\small
$>$ with (plots):

$>$ a:=-0.45; b:=1; c1:=-0.2; c2:=-0.15; c3:=1.01; q:=0.8; x2e:=1; x3e:=0.6; x4e:=0;\\[0.1cm]
$>$ with (stats):

$>$ h:=0.01; epsilon:=0.01;
 n:=100:t:=n+2; x1:= array (0 .. n): x2:= array (0 .. n): x3:= array (0 .. n):x4:= array (0 .. n):x5:= array (0 .. n): x1[0]:=epsilon ; x2[0]:=epsilon+ x2e; x3[0]:=epsilon + x3e; x4[0]:=epsilon + x4e; x5[0]:=epsilon;\\[0.1cm]
$>$ for ~j~ from ~1~  by ~1~ to ~n~  do

$>$ x1[j]:= x1[j-1] + h $\wedge$ q *( x4[j-1]$\wedge$ 2 + a*x1[j-1])/GAMMA(q+1);

 x2[j]:= x2[j-1] + h $\wedge$ q *(2*(x5[j-1]$\wedge$ 2 - x4[j-1]$\wedge$ 2) + c1*(x2[j-1] - x2e))/GAMMA(q+1);

x3[j]:= x3[j-1] + h $\wedge$ q *(-2*x5[j-1]$\wedge$ 2 + c2*(x3[j-1] - x3e))/GAMMA(q+1);

x4[j]:= x4[j-1] + h $\wedge$ q *(x4[j-1]*(x2[j-1]-x1[j-1]) + c3*(x4[j-1] - x4e))/GAMMA(q+1);

x5[j]:= x5[j-1] + h $\wedge$ q *(x5[j-1]*(x3[j-1]-x2[j-1])+ b*x5[j-1])/GAMMA(q+1);\\[0.1cm]
od:

$>$ plot (seq([j,x1[j]], j = 0 .. n), style = point, symbol = point, scaling = UNCONSTRAINED);\\
plot (seq([j,x2[j]], j = 0 .. n), style = point, symbol = point, scaling = UNCONSTRAINED);\\
plot (seq([j,x3[j]], j = 0 .. n), style = point, symbol = point, scaling = UNCONSTRAINED);\\
plot (seq([j,x4[j]], j = 0 .. n), style = point, symbol = point, scaling = UNCONSTRAINED);\\
plot (seq([j,x5[j]], j = 0 .. n), style = point, symbol = point, scaling = UNCONSTRAINED);}
\begin{Rem}
{\rm  Appyling  $ (4.8) $ and  Maple for the numerical simulation of solution of fractional model $~(4.3)~$ for each set of values for parameters $~ a, b, c_{1}, c_{2}, c_{3},~$ given in the Examples 3.1-3.5, it will be found that the results obtained are valid.}\hfill$\Box$
\end{Rem}
{\bf Conclusions.} This paper presents the 3-dimensional fractional-order Toda lattice with two controls, denoted by $ (2.4).$  The fractional model $ (2.4) $ was studied from fractional differential equations theory point of view: asymptotic stability, determining of sufficient conditions on parameters $~a, b, c_{1}, c_{2}, c_{3}~$ to control the chaos in the controlled fractional system associated to $ (2.4) $ and numerical integration of the fractional model $ (3.3).~$ By choosing the right parameters $ a, b, c_{1}, c_{2} $ and $ c_{3} $ in the fractional model  $~(3.3),~$ this work offers a series of chaotic fractional differential systems. \hfill$\Box$\\
{\bf Acknowledgments.} The author has very grateful to be reviewers for their comments and suggestions.

Author's adress\\[-0.5cm]

Mihai Ivan\\[0.1cm]
West University of Timi\c soara. Seminarul de Geometrie \c si Topologie.\\
Teacher Training Department. Timi\c soara, Romania.\\
E-mail: mihai.ivan@e-uvt.ro\\
\end{document}